\documentclass{amsart}

\usepackage{amsmath,amsthm,amssymb,amsfonts}

\theoremstyle{plain}
\newtheorem{theo+}           {Theorem}      [section]
\newtheorem{prop+}  [theo+]  {Proposition}
\newtheorem{lemm+}  [theo+]  {Lemma}
\newtheorem{cor+}  [theo+]  {Corollary}

\theoremstyle{definition}
\newtheorem{exam+}  [theo+]  {Example}
\newtheorem{rema+}  [theo+]  {Remark}

\newenvironment{theorem}{\begin{theo+}}{\end{theo+}}

\newenvironment{lemma}{\begin{lemm+}}{\end{lemm+}}
\newenvironment{corollary}{\begin{cor+}}{\end{cor+}}
\newenvironment{remark}{\begin{rema+}}{\end{rema+}}

\newcommand{\C}{\mathbb C}
\newcommand{\R}{\mathbb R}
\newcommand{\cmp}{\C\setminus\R_{>0}}

\begin{document}
\baselineskip 18pt
\larger[2]
\title[A bilateral series involving basic hypergeometric functions]
{A bilateral series involving\\ basic hypergeometric functions}
\author{Hjalmar Rosengren}
\address
{Department of Mathematics\\ Chalmers University of Technology and G\"oteborg
 University\\SE-412~96 G\"oteborg, Sweden}
\email{hjalmar@math.chalmers.se}

\dedicatory{\large Dedicated to Mizan Rahman}

\begin{abstract}
We prove a summation formula for a bilateral series whose  terms are
 products of two basic hypergeometric functions. In special cases,
series of this type arise as matrix elements of quantum group representations.

\end{abstract}

\maketitle        

\section{Introduction}  

The object of the present paper is to study bilateral series of a type
that first appeared in the work of Koelink and Stokman \cite{ks}, in
the context of
harmonic analysis on the $\mathrm{SU}(1,1)$ quantum group. They needed 
 to compute sums of the form
\begin{equation}\label{os}\sum_{n=-\infty}^\infty  {}_{2}\phi_1\left(
  {{a,b}
\atop{c}};q,
xq^n\right) {}_{2}\phi_1\left(  {{d,e}
\atop{f}};q,yq^n\right)t^n,\end{equation}
which appeared as matrix elements of quantum $\mathrm{SU}(1,1)$
representations.
Several special cases were considered,
 with different, very technical,
proofs (some of these proofs were found by Mizan Rahman, see the
appendix to \cite{ks}). 

In joint work with Koelink \cite{kr}, we gave an
extension of the summation formulas from \cite{ks}, with a
unified proof. Namely, we showed that under natural conditions of
convergence and the assumption
\begin{equation}\label{cb}abde=cf,\qquad fx=dey\end{equation}
(and thus also $cy=abx$),
the sum in \eqref{os}  can be expressed  as the sum of two ${}_8W_7$
series, or, alternatively, as the sum of three balanced ${}_4\phi_3$ series.
From the quantum algebraic viewpoint, which was only briefly
mentioned in \cite{kr}, this 
computes a general class of matrix elements for the strange, the
complementary and the principal unitary series of quantum  $\mathrm{SU}(1,1)$.
The proof in \cite{kr} is similar to, though not an extension of, one of
Rahman's proofs in \cite{ks}. 
In particular, it is quite technical and involves 
 rather non-obvious applications of 
$q$-series transformations. 

In a related paper \cite{s}, Stokman considered
 matrix elements for the principal unitary
series  and gave a  very simple proof
for their expression as ${}_8W_7$ series. The idea is simple
but powerful: the
representations are realized on $\mathrm{L}^2(\mathbb{T})$, where one
wants to compute the scalar product of two particular functions. 
Expanding both functions as Fourier series gives a special case of 
the sum \eqref{os}.
Stokman's observation, which is our starting point,
 is that the integral defining the scalar
product can be computed very easily using residue calculus. 

Our  aim is  to use Stokman's method not only to give a simple proof
 of the summation formula in \cite{kr}, but also to give a far-reaching
 extension of this identity.

{\it Dedication:} 
This paper 
would not have been written without the  pioneering contributions
of Mizan Rahman. 
I'm happy to know him, not only through his work as a master
of identities, but also as a genuinely kind and helpful person.
It is a pleasure to dedicate this paper to him as a small token of
my appreciation.

\section{A bilateral summation}

We will now state the main result of the paper. Throughout, we use the
standard notation of \cite{gr}. The base $q$ is assumed to satisfy $0<|q|<1$.
We denote by ${}_{r+1}\phi_r$ not only
the convergent basic hypergeometric series, but also its
analytic continuation to $\mathbb C\setminus\mathbb R_{\geq 1}$; cf.\
\cite[\S 4.5]{gr}. 

\begin{theorem}\label{t}
Let  $x,\,y\in \cmp$, and let the other parameters below satisfy
\begin{gather}\label{acl}\max_{ij}(|a_ic_j|)<|t|<1,\\
\label{acr}|qyb_1\dotsm b_{k}|<|xa_1\dotsm a_{k+1}|.
\end{gather}
Then the following identity holds:
\begin{multline}\label{gs}\sum_{n=-\infty}^\infty  {}_{k+1}\phi_k\left(
  {{a_1,\dots,a_{k+1}}
\atop{b_1,\dots,b_k}};q,
xq^n\right) {}_{l+1}\phi_l\left(  {{c_1,\dots,c_{l+1}}
\atop{d_1,\dots,d_l}};q,yq^n\right)t^n\\
\begin{split}&=\frac{(q,xt,q/xt,a_1,\dots,a_{k+1},b_1/t,\dots,b_k/t;q)_\infty}
{(t,x,q/x,a_1/t,\dots,a_{k+1}/t,b_1,\dots,b_k;q)_\infty}\\
&\qquad\times{}_{k+l+2}\phi_{k+l+1}\left({{t,qt/b_1,\dots,qt/b_k,
    c_1,\dots,c_{l+1}}\atop{qt/a_1,\dots,qt/a_{k+1},
    d_1,\dots,d_{l}}};q,\frac{qyb_1\dotsm b_k}{xa_1\dotsm a_{k+1}}\right)\\
&\quad+\frac{(q,a_1x,q/a_1x,yt/a_1,qa_1/yt,a_2,\dots,a_{k+1},b_1/a_1,\dots,
    b_k/a_1;q)_\infty}
{(x,q/x,y,q/y,t/a_1,a_2/a_1,\dots,a_{k+1}/a_1,b_1,\dots,b_k;q)_\infty}\\
&\qquad\times\frac{(c_1,\dots,c_{l+1},a_1d_1/t,\dots, a_1d_l/t;q)_\infty}{(
    a_1c_1/t,\dots, a_1c_{l+1}/t,d_1,\dots,d_l;q)_\infty}\\
&\qquad\times {}_{k+l+2}\phi_{k+l+1}\left({{a_1,qa_1/b_1,\dots,qa_1/b_k,
    a_1c_1/t,\dots,a_1c_{l+1}/t}\atop{qa_1/t,qa_1/a_2,\dots,qa_1/a_{k+1},
    a_1d_1/t,\dots,a_1d_{l}/t}}; 
q,\frac{qyb_1\dotsm b_k}{xa_1\dotsm a_{k+1}}\right)\\
&\quad+\operatorname{idem}(a_1;a_2,\dots,a_{k+1}).
 \end{split}\end{multline}
\end{theorem}

Here, 
$\operatorname{idem}(a_1;a_2,\dots,a_{k+1})$ denotes the sum of
the $k$ terms obtained by interchanging $a_1$ with each of 
$a_2,\dots, a_{k+1}$. As is customary, we   implicitly assume
 that the parameters are such
that one never divides by zero.
Note also that interchanging the roles of ${}_{k+1}\phi_k$ and
${}_{l+1}\phi_l$ gives an alternative expression for the series, which is
 valid when \eqref{acr} is
replaced by the condition
$|qxd_1\dotsm d_{l}|<|yc_1\dotsm c_{l+1}|$. These two expressions are
related by \cite[(4.10.10)]{gr}.

In \cite{kr}, we computed the sum 
when $k=l=1$ and, crucially to the methods used
there,  \eqref{cb} holds. The latter condition 
means exactly that the three ${}_4\phi_3$ series on the right-hand
side of \eqref{gs} are balanced. In this case Theorem \ref{t} reduces
to \cite[(3.13)]{kr}, rather than the alternative expression as a sum
of two ${}_8W_7$ series given in \cite[Proposition 3.1]{kr}.

Our main tool is the following lemma from \cite{kr}, where it was used
only in the case $k=5$. Here we need the general case. 

\begin{lemma}\label{l}
For 
\begin{equation}\label{a}\max_i(|a_i|)<|t|<1\end{equation}
 and  $x\in\cmp$ one has
\begin{multline}\label{li}
\sum_{n=-\infty}^\infty
 {}_{k+1}\phi_k\left(
{{a_1,\ldots,a_{k+1}}\atop{b_1,\ldots,b_k}};q,xq^n\right) t^n\\
=\frac{(q,a_1,\ldots,a_{k+1},xt,q/xt,b_1/t,\ldots,b_k/t;q)_\infty}
{(x,q/x,b_1,\ldots,b_k,t,a_1/t,\ldots,a_{k+1}/t;q)_\infty}.
\end{multline}
\end{lemma}

In \cite{kr}, this was proved using the explicit formula for 
analytic  continuation of ${}_{k+1}\phi_k$ series. It is possible to
give a much simpler proof using integral representations. Namely, the
coefficient of $t^n$ in the Laurent expansion of the right-hand side 
in the annulus \eqref{a} is given by
\begin{equation}\label{lin}\frac{(q,a_1,\ldots,a_{k+1};q)_\infty}
{(x,q/x,b_1,\ldots,b_k;q)_\infty}\frac{1}{2\pi i}
\int\frac{(xt,q/xt,b_1/t,\ldots,b_k/t;q)_\infty}
{(t,a_1/t,\ldots,a_{k+1}/t;q)_\infty}\,t^{-n-1}\,dt,\end{equation}
where the integral is over a positively oriented contour encircling
the origin inside the annulus. This is an integral of the form 
\cite[(4.9.4)]{gr}, which is computed there using residue calculus.
 If $|xq^n|<1$, 
its value is given by  \cite[(4.10.9)]{gr}
as the ${}_{k+1}\phi_k$ series in \eqref{li}. Since \eqref{lin} is
analytic in $x$, this also holds for $|xq^n|\geq 1$ in the sense of
analytic continuation. 

\begin{remark}\label{r}
More generally, \cite[(4.10.9)]{gr} may be used to express the Laurent
coefficients  of
\begin{equation}\label{gp}
\prod_{i=1}^m\frac{(\alpha_i t;q)_\infty}{(\gamma_i t;q)_\infty}
\prod_{i=1}^n\frac{(\beta_i/ t;q)_\infty}{(\delta_i/ t;q)_\infty}
\end{equation}
in the annulus 
\begin{equation}\label{ga}\max_i|\delta_i|<|t|<\min_i(1/|\gamma_i|)
\end{equation}
as sums of analytically continued basic hypergeometric series.
\end{remark}

\begin{proof}[Proof of \emph{Theorem \ref{t}}]
 Similarly as in \cite{kr}, 
 it is easy to check that \eqref{acl}
is the natural condition for
absolute convergence of the left-hand side of \eqref{gs}. 
We first rewrite this series
as an integral. For this we  make the preliminary assumption 
that $t$ is real with
\begin{equation}\label{cc}\max_{ij}(|a_i|,|c_{j}|)<t^{1/2}<1.\end{equation} 
Writing
$f_k(t;a,b,x)$ for either side of \eqref{li}, we  consider
the integral
$$\frac{1}{2\pi i}\int_{|z|=t^{1/2}}f_k(z;a,b,x)\,f_l(\bar z;c,d,y)\,
\frac{dz}{z}.$$
Using the expression on the 
left-hand side of \eqref{li} together with orthogonality of the
monomials, we find that it equals the sum in  \eqref{gs}.

To compute the integral,  we plug in the expressions from the right-hand
side of \eqref{li}. The integrand is then of the form
\eqref{gp}, with
$\alpha=(x,q/yt,d_1/t,\dots,d_l/t)$, $\beta=(q/x,yt,b_1,\dots,b_k)$, 
$\gamma=(1,c_1/t,\dots,c_{l+1}/t)$, $\delta=(t,a_1,\dots,a_{k+1})$.
Note that \eqref{ga} and 
 \eqref{cc} are equivalent. Thus, we again have 
 an integral of the form  \cite[(4.9.4)]{gr}.
The condition \eqref{acr} is  precisely
\cite[(4.10.2)]{gr}, which 
ensures that the singularity at $z=0$ does not contribute
to the integral. The value of the integral is then given by    
\cite[(4.10.8)]{gr} as a sum of $k+2$ terms, each being a  
 ${}_{k+l+4}\phi_{k+l+3}$ series. However, because of the  condition
$\alpha_1\beta_1=\alpha_2\beta_2=q$, they all  reduce to 
type ${}_{k+l+2}\phi_{k+l+1}$. 
This proves
 Theorem \ref{t} under the assumption \eqref{cc}.
By analytic continuation in $t$,
this may be replaced with the weaker condition
 \eqref{acl}.
\end{proof}

\begin{remark}
Using instead of Lemma \ref{l} the Laurent expansion of the general
product \eqref{gp} gives a generalization of Theorem \ref{t}, where
the two basic hypergeometric series on the left-hand side of
\eqref{gs} are replaced by finite sums of such series. 
\end{remark}

In the case $l=0$, one may use the $q$-binomial theorem to sum the
series ${}_1\phi_0$ in Theorem \ref{t}.  The condition $y\in \cmp$ is
then superfluous. We find this special case interesting enough
to write out explicitly. Compared to Theorem \ref{t}  we have made the change 
of  variables $y\mapsto c$, $c_1\mapsto d/c$.

\begin{corollary}\label{c}
Let  $x\in \cmp$, and assume that
$$\max_{i}(|a_id/c|)<|t|<1,\qquad
|qcb_1\dotsm b_{k}|<|xa_1\dotsm a_{k+1}|.$$
Then the following identity holds:
\begin{multline*}\sum_{n=-\infty}^\infty \frac{(c)_n}{(d)_n} \,
{}_{k+1}\phi_k\left({{a_1,\dots,a_{k+1}}\atop{b_1,\dots,b_k}};q,
xq^n\right)t^n\\
\begin{split}&=\frac{(q,xt,q/xt,c,
a_1,\dots,a_{k+1},b_1/t,\dots,b_k/t;q)_\infty}
{(t,x,q/x,d,a_1/t,\dots,a_{k+1}/t,b_1,\dots,b_k;q)_\infty}\\
&\qquad\times{}_{k+2}\phi_{k+1}\left({{t,d/c,qt/b_1,\dots,qt/b_k
    }\atop{qt/a_1,\dots,qt/a_{k+1}
    }};q,\frac{qcb_1\dotsm b_k}{xa_1\dotsm a_{k+1}}\right)\\
&\quad+\frac{(q,d/c,a_1x,q/a_1x,ct/a_1,qa_1/ct,a_2,\dots,a_{k+1},b_1/a_1,\dots,
    b_k/a_1;q)_\infty}{(x,q/x,q/c,d,a_1d/ct,t/a_1,a_2/a_1,\dots,
a_{k+1}/a_1,b_1,\dots,b_k;q)_\infty}\\
&\qquad\times {}_{k+2}\phi_{k+1}\left({{a_1,a_1d/ct,qa_1/b_1,\dots,qa_1/b_k
   }\atop{qa_1/t,qa_1/a_2,\dots,qa_1/a_{k+1}}}; 
q,\frac{qcb_1\dotsm b_k}{xa_1\dotsm a_{k+1}}\right)\\
&\quad+\operatorname{idem}(a_1;a_2,\dots,a_{k+1}).
 \end{split}\end{multline*}
\end{corollary}

Note that in the case $c=d$, the first ${}_{k+2}\phi_{k+1}$ on the
right-hand side  reduces to $1$ and the remaining terms to $0$. 
Thus, we recover Lemma \ref{l}. We also remark 
that if we  choose $k=0$ in Corollary \ref{c} and  replace $x\mapsto a$,
$a_1\mapsto b/a$, we obtain the transformation formula
\begin{equation*}\begin{split}{}_2\psi_2\left({{a,c}\atop{b,d}}\,;q,t\right)&=
\frac{(q,b/a,c,at,q/at;q)_\infty}{(q/a,b,d,t,b/at;q)_\infty}\,
{}_2\phi_1\left({{t,d/c}\atop{aqt/b}};q,\frac{qc}{b}\right)\\
&\quad + \frac{(q,q/b,d/c,act/b,qb/act;q)_\infty}
{(q/a,q/c,d,at/b,bd/act;q)_\infty}\,
{}_2\phi_1\left({{b/a,bd/act}\atop{qb/at}};q,\frac{qc}{b}\right).
\end{split}\end{equation*}
This identity is also obtained by choosing $r=2$, $c_1=qa_2$, $c_2=b_2$ in 
\cite[(5.4.3)]{gr}, and then applying \cite[(III.1)]{gr} to both
${}_{2}\phi_1$ series on the right-hand side.

\end{document}